\documentclass[11pt]{article}
\usepackage{amsfonts,amscd,amssymb, amsmath}
\newtheorem{definition}{Definition}[section]

\newtheorem{proposition}[definition]{Proposition}

\newtheorem{remark}[definition]{Remark}

\newtheorem{theorem}[definition]{Theorem}

\def\rawo\lonra{\longrightarrow}

\allowdisplaybreaks[4]

\newenvironment{proof}{{\it Proof.}}{\hfill $ \square $ \vskip 4mm}
\begin{document}
\title{Alternative twisted tensor products and Cayley algebras}
\author{Helena Albuquerque\thanks{The first author was partially supported 
by the Centre for Mathematics of the University of Coimbra (CMUC)} \\
Departamento de Matem\'{a}tica, Universidade de Coimbra, \\
3001-454 Coimbra, Portugal\\
%Departamento de Matematica-Faculdade de Ciencias e Technologia,\\
%Universidade de Coimbra, Apartado 3008, 3000 Coimbra, Portugal\\
e-mail: lena@mat.uc.pt 
\and Florin Panaite\thanks {Research
started while the second author was visiting the Centre for 
Mathematics of   
the University of Coimbra (CMUC) supported by a postdoctoral fellowship  
offered by CMUC, and finished while this author was a member of the 
CNCSIS project ``Hopf algebras, cyclic homology and monoidal categories'', 
contract nr. 560/2009, CNCSIS code $ID_{-}69$.}\\
Institute of Mathematics of the
Romanian Academy\\
PO-Box 1-764, RO-014700 Bucharest, Romania\\
e-mail: Florin.Panaite@imar.ro}
\date{}
\maketitle

\begin{abstract}
We introduce what we call {\em alternative twisted tensor products} for   
not necessarily associative algebras, as a common generalization of 
several different constructions: the  
{\em Cayley-Dickson process}, the {\em Clifford process} and the twisted 
tensor product of two associative algebras, one of them being commutative. 
We show that some very basic facts concerning the Cayley-Dickson process  
(the equivalence between the two different formulations of it and the 
lifting of the involution)   
are particular cases of general results about alternative twisted tensor 
products of algebras. As a class of examples of alternative twisted tensor 
products, we introduce a {\em tripling process} for an algebra endowed 
with a strong involution, containing the Cayley-Dickson doubling as a 
subalgebra and sharing some of its basic properties.
\end{abstract}
%%%%%%%%%%%%%%%%%%%%%%%%%%%%%%%%%%%%%%%%%%%%%%%%%%%
\section*{Introduction}
%%%%%%%%%%%%%%%%%%%%%%%%%%%%%%%%%%%%%%%%%%%%%%%%%%%%
${\;\;\;\;}$
If $A$ and $B$ are associative algebras and 
$R:B\otimes A\rightarrow A\otimes B$ is a linear map  
satisfying a certain list of axioms, then 
$A\otimes B$ becomes an associative algebra with a   
multiplication defined in terms of $R$ and the multiplications of $A$ and 
$B$. This  
construction appeared in several contexts and under different names.  
Following \cite{cap} we call such an $R$ a {\it twisting map} and the 
algebra   
structure on $A\otimes B$ afforded by it the {\it twisted tensor product} 
of $A$ and $B$ and denote it by $A\otimes _RB$ (if $R$ is the usual flip 
map then $A\otimes _RB$ coincides with the usual tensor product of 
algebras). The twisted tensor product of associative 
algebras can be regarded as a representative for the cartesian product 
of noncommutative spaces, better suited than the ordinary tensor product, 
see \cite{cap}, \cite{jlpvo}, \cite{lpvo} for a detailed discussion and 
references. Examples of twisted tensor products are provided by classical 
ring theory (crossed products with trivial cocycles, 
Ore extensions with trivial derivations), Hopf algebra theory (smash  
products, diagonal crossed products as in \cite{bpvo}, \cite{hn}), 
noncommutative geometry etc.

An important feature of this construction is that it allows to build new 
algebra structures starting with two given algebras. From this point of view, 
and also having in mind the emerging 
{\em nonassociative geometry} (cf. \cite{akm}, \cite{bm}, 
\cite{majidgauge}), regarded as  
a further extension of noncommutative geometry, with the 
``coordinate algebra'' allowed to be nonassociative, it is natural to try 
to consider analogues of twisted tensor products for more general 
classes of algebras than the associative ones. It turns out that it is 
straightforward to do this for the so-called {\em quasialgebras} as 
named in \cite{AM1}, that is algebras in monoidal categories (since the 
associativity constraints are allowed to be nontrivial, the class of 
quasialgebras contains various classes of nonassociative algebras, for 
instance the octonions and the other Cayley algebras, cf. \cite{AM1}). 
Twisted tensor products of certain quasialgebras have been studied by the 
authors in \cite{AP}. 

In this paper we will introduce a different kind of twisted tensor 
product for nonassociative algebras, having as motivating example the 
Cayley-Dickson process. We recall that if 
$B$ is an algebra endowed with an involution 
$\sigma :B\rightarrow B$ and $q$ is a nonzero element in the base field 
$K$, the Cayley-Dickson algebra $\overline{B}(q)$ is an algebra structure 
on $B\oplus B$, whose elements are written uniquely as   
$a+vb$, with $a, b\in B$, where $v$ is a notational device with 
$v^2=q1$, and whose multiplication is defined by   
\begin{eqnarray*}
&&(a+vb)(c+vd)=(ac+qd\sigma (b))+v(\sigma (a)d+cb), \;\;\;\forall \;\;
a, b, c, d\in B. 
\end{eqnarray*}
If we denote by $C(K, q)$ the 2-dimensional algebra $C(K, q):=K[v]/(v^2-q)$,  
it is clear that $\overline{B}(q)$ may be identified as vector space with 
$C(K, q)\otimes B$, and the multiplication of $\overline{B}(q)$ depends 
somehow on the multiplications of $C(K, q)$ and $B$. Our original purpose  
was to make exact sense of this claim; more precisely, we wanted to 
express $\overline{B}(q)$ as some sort of twisted tensor product between 
$C(K, q)$ and $B$.

The general construction that we introduce, with $\overline{B}(q)$ as the   
motivating example, is called {\em alternative twisted tensor product} 
between two nonassociative algebras $A$ and $B$ and is defined by means 
of a so-called {\em alternative twisting map} 
$R:B\otimes A\rightarrow A\otimes B$, which is a linear map satisfying 
some axioms resembling the ones for a twisting map. We emphasize that a 
twisting map between two associative algebras $A$ and $B$ in general is 
{\em not} an alternative twisting map. However, if $B$ is moreover 
commutative, the two concepts coincide. 

Our main results (Theorems \ref{main} and \ref{ext}) show that some basic 
features of the Cayley-Dickson process (the fact that it admits two 
different but equivalent descriptions and the fact that the involution 
$\sigma $ can be lifted to an involution on $\overline{B}(q)$), which 
apparently are very peculiar to it, are actually particular cases of 
general results about alternative twisted tensor products of algebras. This 
kind of phenomenon appeared also in \cite{jlpvo}, where it was shown that 
various results from Hopf algebra theory are actually particular cases 
of general results about twisted tensor products of associative algebras. 

As a class of examples of alternative twisted tensor products, we 
introduce a sort of {\em tripling process}: if $B$ is an algebra endowed 
with a strong involution $\sigma $ and $q, r$ are nonzero scalars, we  
construct an algebra $\overline{B}(q, r)$, whose dimension is 
$3\cdot dim (B)$, containing the Cayley-Dickson algebras $\overline{B}(q)$ 
and $\overline{B}(r)$ as subalgebras. Unlike $\overline{B}(q)$ (which is 
alternative when $B$ is associative), the algebras 
$\overline{B}(q, r)$ are {\em never} alternative algebras. But, 
exactly as for Cayley algebras, $\overline{B}(q, r)$ is always 
power-associative, it is flexible if and only if $B$ is flexible and if the 
norm on $B$ is nondegenerate then so is the norm on 
$\overline{B}(q, r)$.       
%%%%%%%%%%%%%%%%%%%%%%%%%%%%%%%%%%%%%%%%%%%%%%%%%%%%%%%%%%%%%
\section{Alternative twisted tensor products of algebras}
%%%%%%%%%%%%%%%%%%%%%%%%%%%%%%%%%%%%%%%%%%%%%%%%%%%%%%%%%%%%
\setcounter{equation}{0}
%%%%%%%%%%%%%%%%%%%%%%%%%%%%%%%%%%%%%%%%%%%%%%%%%
${\;\;\;\;}$
In what follows, by ``algebra'' we will mean a  
not necessarily associative  
algebra with unit over a field $K$ (of characteristic $\neq 2$ if 
necessary); all algebra maps are assumed to be unital.  
For an algebra $A$ we will denote its multiplication by 
$\mu _A$ and its unit by $1_A$. All algebras that will appear will be 
considered together with a given and fixed decomposition as a direct sum 
of vector spaces $A=K\cdot 1_A\oplus A_0$ (our constructions will depend 
on the decomposition). If $q\in K$, $q\neq 0$, the decomposition for 
the algebra $C(K, q):=K[v]/(v^2-q)$ will always be 
$C(K, q)=K\cdot 1\oplus K\cdot v$.  
If $V$, $W$ are vector spaces, we denote by $\tau _{V, W}:V\otimes W 
\rightarrow W\otimes V$ the flip map $v\otimes w\mapsto w\otimes v$.

We recall from \cite{cap}, \cite{vandaele} the construction of the 
twisted tensor product of associative algebras. Let $A$ and $B$ be two 
associative algebras and $R:B\otimes A\rightarrow A\otimes B$ a 
linear map, with Sweedler-type notation $R(b\otimes a)=a_R\otimes b_R$, 
for $a\in A$, $b\in B$. Then $R$ is called 
{\em twisting map} if the following conditions are satisfied (we denote by 
$r$ another copy of $R$): 
\begin{eqnarray*}
&&R(1_B\otimes a)=a\otimes 1_B ,\;\;\;
R(b\otimes 1_A)=1_A\otimes b, \\
&&R(b\otimes aa')=a_Ra'_r\otimes (b_R)_r,\\ 
&&R(bb'\otimes a)=(a_R)_r\otimes b_rb'_R, 
\end{eqnarray*} 
for all $a, a'\in A$ and $b, b'\in B$. If we introduce a multiplication 
on $A\otimes B$, by 
\begin{eqnarray*}
&&(a\otimes b)(a'\otimes b')=aa'_R\otimes b_Rb', \;\;\;\forall \;\;
a, a'\in A, b, b'\in B, 
\end{eqnarray*}
then this multiplication is associative with unit $1_A\otimes 1_B$. 
This algebra structure on $A\otimes B$ is called the {\em twisted tensor 
product} afforded by $R$ and is denoted by $A\otimes _RB$.

Let $A$ and $B$ be two algebras and $R:B\otimes A\rightarrow A\otimes B$ a 
linear map, with Sweedler-type notation $R(b\otimes a)=a_R\otimes b_R$, 
for $a\in A$, $b\in B$. Assume that the following conditions are satisfied 
(we denote by $r$ another copy of $R$):  
\begin{eqnarray}
&&R(1_B\otimes a)=a\otimes 1_B ,\;\;\;
R(b\otimes 1_A)=1_A\otimes b, \;\;\;\forall \;\;a\in A, \;b\in B, 
\label{atm1} \\
&&R(b\otimes aa')=a_Ra'_r\otimes (b_R)_r,\;\;\;\forall \;\;
a,a'\in A, \; b\in B, \label{atm2} \\ 
&&R(bb'\otimes a)=(a_R)_r\otimes b'_Rb_r, \;\;\;\forall \;\;
a\in A_0, \;b, b'\in B. \label{atm3}
\end{eqnarray} 

We consider the multiplication on the vector space $A\otimes B$ 
uniquely defined by the formulae 
\begin{eqnarray*}
&&(1_A\otimes b)(a'\otimes b')=a'_R\otimes b_Rb', \;\;\;\forall \;\;
a'\in A, \; b, b'\in B, \\
&&(a\otimes b)(a'\otimes b')=aa'_R\otimes b'b_R, \;\;\;\forall \;\;
a\in A_0, \; a'\in A, \;b, b'\in B.
\end{eqnarray*}

This algebra structure on $A\otimes B$ will be called an 
{\bf alternative twisted tensor product} and will be denoted by 
$A\overline{\otimes }_RB$; the map $R$ satisfying the conditions 
(\ref{atm1})--(\ref{atm3}) will be called an {\bf alternative twisting 
map}. Clearly $1_A\otimes 1_B$ is the unit for  
$A\overline{\otimes }_RB$. 
If the algebras $A$ and $B$ are associative and $B$ is commutative 
then the alternative twisted tensor product $A\overline{\otimes }_RB$  
coincides with the usual twisted tensor product $A\otimes _RB$ and so  
is an associative algebra.  
\begin{remark}
If $A\overline{\otimes }_RB$ is an alternative twisted tensor product 
and $B$ is commutative, then the multiplication of 
$A\overline{\otimes }_RB$ does not depend on the decomposition 
$A=K\cdot 1_A\oplus A_0$.  
\end{remark}

We recall now the Cayley-Dickson process, see for instance \cite{schafer}. 
Let $B$ be an algebra and $\sigma :B\rightarrow B$ an involution,  
i.e. an algebra antiautomorphism with $\sigma ^2=id_B$. 
We fix $q\in K$, $q\neq 0$, and define $\overline{B}(q):=B\oplus B$ as  
vector space; we write an element of $\overline{B}(q)$ uniquely as  
$a+vb$, with $a, b\in B$, where, as above, $C(K, q)=K[v]/(v^2-q)$.   
Define a multiplication on $\overline{B}(q)$, by   
\begin{eqnarray*}
&&(a+vb)(c+vd)=(ac+qd\sigma (b))+v(\sigma (a)d+cb), \;\;\;\forall \;\;
a, b, c, d\in B. 
\end{eqnarray*}
This algebra $\overline{B}(q)$ is said to have  
been obtained from $B$ by the Cayley-Dickson process.

Consider the linear map $R:B\otimes C(K, q)\rightarrow 
C(K, q)\otimes B$ uniquely defined by the formulae  
\begin{eqnarray}
&&R(b\otimes 1)=1\otimes b, \;\;\;R(b\otimes v)=v\otimes \sigma (b), \;\;
\forall \;\;b\in B. \label{Rmap}
\end{eqnarray}   

Then one can check that $R$ is an alternative twisting map 
and we have an algebra isomorphism   
\begin{eqnarray*}
&&\overline{B}(q)\simeq C(K, q)\overline{\otimes }_RB, \;\;\;a+vb\mapsto  
1\otimes a+v\otimes b, \;\;\;\forall \;\;a, b\in B.
\end{eqnarray*}

Thus, any algebra obtained by applying the 
Cayley-Dickson process may be regarded as an alternative twisted tensor 
product of algebras. In particular, this is the case of the algebra of 
octonions and, by \cite{aepi}, of any division alternative quasialgebra.  

We recall now the so-called {\em Clifford process} as introduced in 
\cite{AM2}. Let $A$ be an algebra and $\sigma :A\rightarrow A$ an 
algebra automorphism which is involutive (i.e. $\sigma ^2=id_A$), let 
$q\in K$, $q\neq 0$ and again $C(K, q)=K[v]/(v^2-q)$. Define   
$Cl(A):=A\oplus A$ as vector space and write an element of $Cl(A)$ 
uniquely as $a+bv$, with $a, b\in A$; then $Cl(A)$ becomes an algebra 
with multiplication given by the formula
\begin{eqnarray*}
&&(a+bv)(c+dv)=(ac+qb\sigma (d))+(ad+b\sigma (c))v, \;\;\;\forall \;\;
a, b, c, d\in A. 
\end{eqnarray*}   
This algebra $Cl(A)$ is said to have been obtained  
from $A$ by the Clifford process. 

Consider the linear map $R:C(K, q)\otimes A\rightarrow A\otimes C(K, q)$ 
uniquely defined by the formulae 
\begin{eqnarray*}
&&R(1\otimes a)=a\otimes 1, \;\;\;R(v\otimes a)=\sigma (a)\otimes v, 
\;\;\;\forall \;\;a\in A.
\end{eqnarray*}
Then one can check that $R$ is an alternative twisting map and we have an 
algebra isomorphism 
\begin{eqnarray*}
&&Cl(A)\simeq A\overline{\otimes }_RC(K, q), \;\;\;a+bv\mapsto 
a\otimes 1+b\otimes v, \;\;\;\forall \;\;a, b\in A.
\end{eqnarray*}

If $A$ is associative, since $C(K, q)$ is associative and commutative 
it follows that $A\overline{\otimes }_RC(K, q)$ is a usual twisted 
tensor product of associative algebras, so $Cl(A)$ is an associative algebra 
(this was noted in \cite{AM2} too). By \cite{AM2}, the usual Clifford 
algebras may be obtained from the field $K$ by iterating the Clifford 
process. Thus, the alternative twisted tensor product of algebras 
as introduced above provides a common generalization of Clifford algebras 
and Cayley algebras.

Let $A\overline{\otimes }_RB$ be an alternative twisted tensor product of 
algebras. Clearly the maps 
\begin{eqnarray*}
&&A\rightarrow A\overline{\otimes }_RB, \;\;\;a\mapsto a\otimes 1_B, \\
&&B\rightarrow A\overline{\otimes }_RB. \;\;\;b\mapsto 1_A\otimes b, 
\end{eqnarray*}  
are algebra maps, and $(a\otimes 1_B)(1_A\otimes b)=a\otimes b$ in 
$A\overline{\otimes }_RB$, for all $a\in A$, $b\in B$. We have already 
seen that if $A$ is associative and $B$ is associative and commutative 
then $A\overline{\otimes }_RB$ is associative. We can prove a converse of 
this, generalizing the well-known fact that the algebra $\overline{B}(q)$ 
obtained from $B$ by the Cayley-Dickson process is associative if and 
only if $B$ is associative and commutative:
\begin{proposition}  
If $A_0\neq 0$ (i.e. $dim (A)\geq 2$) and $A\overline{\otimes }_RB$ is 
associative, then $A$ is associative and $B$ is associative and commutative. 
Thus, in this case $A\overline{\otimes }_RB$ is a usual twisted tensor 
product of associative algebras.
\end{proposition}
\begin{proof}
Obviously $A$ and $B$ are associative, as they are embedded as algebras into 
$A\overline{\otimes }_RB$. Let $a\in A_0$, $a\neq 0$ and $b, b'\in B$. 
We have: 
\begin{eqnarray*}
&&(a\otimes 1_B)[(1_A\otimes b)(1_A\otimes b')]=a\otimes bb', \\
&&[(a\otimes 1_B)(1_A\otimes b)](1_A\otimes b')=(a\otimes b)(1_A\otimes b')=
a\otimes b'b,
\end{eqnarray*}
and so the associativity of $A\overline{\otimes }_RB$ implies $bb'=b'b$.
\end{proof}
%%%%%%%%%%%%%%%%%%%%%%%%%%%%%%%%%%%%%%%%%%%
\section{An isomorphism theorem}
%%%%%%%%%%%%%%%%%%%%%%%%%%%%%%%%%%%%%%%%%%%%
\setcounter{equation}{0}
%%%%%%%%%%%%%%%%%%%%%%%%%%%%%%%%%%%%%%%%%%%%%%%%%
${\;\;\;\;}$ 
It is known (see for instance \cite{wene}) that the Cayley-Dickson 
process admits a second description, different from the one presented above, 
but equivalent to it, which we recall now. 

We begin with  
an algebra $B$ and an involution $\sigma $ on it, we fix $q\in K$,  
$q\neq 0$ and consider again the algebra $C(K, q)=K[v]/(v^2-q)$. 
Define $\underline{B}(q):= 
B\oplus B$ as vector space, and write an element of $\underline{B}(q)$  
uniquely as $a+bv$, with $a, b\in B$. Define a multiplication on 
$\underline{B}(q)$ by 
\begin{eqnarray*}
&&(a+bv)(c+dv)=(ac+q\sigma (d)b)+(b\sigma (c)+da)v,  \;\;\;\forall \;\;
a, b, c, d\in B.  
\end{eqnarray*}    
Then one can see that we have an algebra isomorphism 
\begin{eqnarray}
&&\underline{B}(q)\simeq \overline{B}(q), 
\;\;\;a+bv \mapsto a+v\sigma (b), \;\;\;
\forall \;\;a, b\in B. \label{isocd}
\end{eqnarray}  

Our next aim is to give an interpretation of this fact in terms of 
alternative twisted tensor products of algebras. To begin with, we note that 
if $A\overline{\otimes }_RB$ is such an alternative twisted tensor 
product, the algebras $A$ and $B$ do {\em not} play a symmetric r\^{o}le 
in its construction. Thus, we construct first a new kind of 
product, a sort of mirror image of $A\overline{\otimes }_RB$. Namely, we 
begin with two algebras $C$, $D$ and a linear map 
$P:D\otimes C\rightarrow C\otimes D$, with Sweedler-type notation as before, 
and we assume that the following conditions are satisfied (denote by $p$ 
another copy of $P$):    
\begin{eqnarray}
&&P(1_D\otimes c)=c\otimes 1_D ,\;\;\;
P(d\otimes 1_C)=1_C\otimes d, \;\;\;\forall \;\;c\in C, \;d\in D, 
\label{aatm1} \\
&&P(d\otimes cc')=c'_pc_P\otimes (d_P)_p,\;\;\;\forall \;\;
c,c'\in C, \; d\in D_0, \label{aatm2} \\ 
&&P(dd'\otimes c)=(c_P)_p\otimes d_pd'_P, \;\;\;\forall \;\;
c\in C, \;d, d'\in D. \label{aatm3}
\end{eqnarray}  

We consider the multiplication on the vector space $C\otimes D$  
uniquely defined by the formulae 
\begin{eqnarray*}
&&(c\otimes d)(c'\otimes 1_D)=cc'_P\otimes d_P, \;\;\;\forall \;\;
c, c'\in C, \; d\in D, \\
&&(c\otimes d)(c'\otimes d')=c'_Pc\otimes d_Pd', \;\;\;\forall \;\;
c, c'\in C, \; d\in D, \;d'\in D_0.
\end{eqnarray*}

This algebra structure on $C\otimes D$, whose unit is $1_C\otimes 1_D$, 
will be denoted by $C\underline{\otimes }_PD$.   
If the algebras $C$ and $D$ are associative and $C$ is commutative  
then $C\underline{\otimes }_PD$   
coincides with the usual twisted tensor product $C\otimes _PD$ and so   
is an associative algebra.  

Consider now again an algebra $B$ with an involution $\sigma $ and  
the linear map  
$P:C(K, q)\otimes B\rightarrow B\otimes C(K, q)$ uniquely 
defined by the formulae 
\begin{eqnarray}
&&P(1\otimes b)=b\otimes 1, \;\;\;P(v\otimes b)=\sigma (b)\otimes v, \;\;
\forall \;\;b\in B. \label{Pmap}
\end{eqnarray}  
Then $P$ satisfies the conditions   
(\ref{aatm1})--(\ref{aatm3}) and we have an algebra isomorphism 
\begin{eqnarray*}
&&\underline{B}(q)\simeq B\underline{\otimes }_PC(K, q), \;\;\;a+bv\mapsto  
a\otimes 1+b\otimes v, \;\;\;\forall \;\;a, b\in B.
\end{eqnarray*}

We can regard now the isomorphism (\ref{isocd}) between the two 
formulations $\underline{B}(q)$ and $\overline{B}(q)$ of the Cayley-Dickson  
process as follows: we have an algebra isomorphism 
$\varphi :B\underline{\otimes }_PC(K, q)\simeq C(K, q)
\overline{\otimes }_RB$, given by 
\begin{eqnarray}
&&\varphi (b\otimes 1)=1\otimes b, \;\;\;\varphi (b\otimes v)=
v\otimes \sigma (b), \;\;\;\forall \;\;b\in B, \label{izom}
\end{eqnarray} 
where $R$ is the map given by (\ref{Rmap}) and $P$ is the map given by 
(\ref{Pmap}). 

It turns out that this result is a particular case of a general property   
of alternative twisted tensor products of algebras:  
\begin{theorem}\label{main}
Let $A\overline{\otimes }_RB$ be an alternative twisted tensor product 
of algebras such that the map $R$ is bijective, with inverse denoted 
by $P:A\otimes B\rightarrow B\otimes A$. Assume that the following 
conditions are satisfied: 
\begin{eqnarray}
&&(id_A\otimes \tau _{B, B})\circ (R\otimes id_B)\circ (id_B\otimes R)
=(R\otimes id_B)\circ (id_B\otimes R)\circ (\tau _{B, B}\otimes id_A), 
\label{braid} \\
&&R(B\otimes A_0)=A_0\otimes B.  \label{cucu}
\end{eqnarray}
Then the map $P$ satisfies the conditions (\ref{aatm1})--(\ref{aatm3}) 
(for $C=B$ and $D=A$) and $R:B\underline{\otimes }_PA\rightarrow 
A\overline{\otimes }_RB$ is an algebra isomorphism. 
\end{theorem}
\begin{proof}
It is obvious that $P$ satisfies (\ref{aatm1}) because $R$ satisfies 
(\ref{atm1}), and $P$ satisfies (\ref{aatm3}) because $R$ satisfies  
(\ref{atm2}). We have to prove that $P$ satisfies (\ref{aatm2}). The 
condition (\ref{atm3}) for $R$ may be written on $B\otimes B\otimes A_0$ as 
\begin{eqnarray*}
R\circ (\mu _B\otimes id_A)&=&(id_A\otimes \mu _B)\circ (id_A\otimes 
\tau _{B, B})\circ (R\otimes id_B)\circ (id_B\otimes R), 
\end{eqnarray*}
which by using (\ref{braid}) becomes 
\begin{eqnarray*}
R\circ (\mu _B\otimes id_A)&=&(id_A\otimes \mu _B)\circ 
(R\otimes id_B)\circ (id_B\otimes R)
\circ (\tau _{B, B}\otimes id_A), 
\end{eqnarray*} 
which by composing with the appropriate maps and using (\ref{cucu}) becomes  
\begin{eqnarray*}
P\circ (id_A\otimes \mu _B)&=&(\mu _B\otimes id_A)\circ (\tau _{B, B}
\otimes id_A)\circ (id_B\otimes P)\circ (P\otimes id_B) 
\end{eqnarray*} 
as an equality of maps from $A_0\otimes B\otimes B$ to $B\otimes A_0$, 
and this is exactly (\ref{aatm2}). \\
The only thing left to prove is that 
$R$ is multiplicative. First, it is obvious that 
\begin{eqnarray*}
&&R((b\otimes 1_A)(b'\otimes 1_A))=R(b\otimes 1_A)R(b'\otimes 1_A), \;\;\;
\forall \;\;b, b'\in B.
\end{eqnarray*}
Take now $a\in A_0$, $b, b'\in B$. By making use of (\ref{cucu}), 
we compute:
\begin{eqnarray*}
R((b\otimes a)(b'\otimes 1_A))&=&R(bb'_P\otimes a_P)\\
&\overset{(\ref{atm3})}{=}&((a_P)_R)_r\otimes (b'_P)_Rb_r\\
&\overset{P=R^{-1}}{=}&a_r\otimes b'b_r\\
&=&(a_r\otimes b_r)(1_A\otimes b')\\
&=&R(b\otimes a)R(b'\otimes 1_A). 
\end{eqnarray*}
Take now $a'\in A_0$, $b, b'\in B$. We compute:
\begin{eqnarray*}
R((b\otimes 1_A)(b'\otimes a'))&=&R(b'b\otimes a')\\
&\overset{(\ref{atm3})}{=}&(a'_R)_r\otimes b_Rb'_r\\
&\overset{(\ref{braid})}{=}&(a'_R)_r\otimes b_rb'_R\\
&=&(1_A\otimes b)(a'_R\otimes b'_R)\\
&=&R(b\otimes 1_A)R(b'\otimes a'). 
\end{eqnarray*}
Finally, take $a, a'\in A_0$, $b, b'\in B$. Using again (\ref{cucu}), 
we compute (denote $r$, $\mathcal{R}$ two copies of $R$):
\begin{eqnarray*}
R((b\otimes a)(b'\otimes a'))&=&R(b'_Pb\otimes a_Pa')\\
&\overset{(\ref{atm2})}{=}&(a_P)_Ra'_r\otimes ((b'_Pb)_R)_r\\
&\overset{(\ref{atm3})}{=}&((a_P)_{\mathcal{R}})_Ra'_r\otimes 
(b_{\mathcal{R}}(b'_P)_R)_r\\
&\overset{(\ref{braid})}{=}&((a_P)_{\mathcal{R}})_Ra'_r\otimes 
(b_R(b'_P)_{\mathcal{R}})_r\\
&\overset{P=R^{-1}}{=}&a_Ra'_r\otimes (b_Rb')_r\\
&\overset{(\ref{atm3})}{=}&a_R(a'_r)_{\mathcal{R}}\otimes 
b'_r(b_R)_{\mathcal{R}}\\
&=&(a_R\otimes b_R)(a'_r\otimes b'_r)\\
&=&R(b\otimes a)R(b'\otimes a'),
\end{eqnarray*}
showing that indeed $R$ is multiplicative.
\end{proof}

If we consider again an algebra $B$ with an involution $\sigma $, 
it is easy to see that the map $R$ given by (\ref{Rmap}) satisfies the 
hypotheses of Theorem \ref{main}, and the isomorphism 
$B\underline{\otimes }_PC(K, q)\simeq   
C(K, q)\overline{\otimes }_RB$ provided by the Theorem coincides with  
the isomorphism $\varphi $ given by (\ref{izom}). 

Let us also note that the condition (\ref{braid}) is a particular case of 
the hexagon (or braid) relation which arose in \cite{jlpvo} in the 
context of iterated twisted tensor products of associative algebras.    
%%%%%%%%%%%%%%%%%%%%%%%%%%%%%%%%%%%%%%%%%%%
\section{Lifting involutions to alternative twisted tensor products}
%%%%%%%%%%%%%%%%%%%%%%%%%%%%%%%%%%%%%%%%%%%%
\setcounter{equation}{0}
%%%%%%%%%%%%%%%%%%%%%%%%%%%%%%%%%%%%%%%%%%%%%%%%%
${\;\;\;\;}$
Let $B$ be an algebra and $\sigma :B\rightarrow B$ an 
involution. It is well-known (see \cite{schafer}) that the map 
\begin{eqnarray}
&&\overline{\sigma }:\overline{B}(q)\rightarrow \overline{B}(q), \;\;\;
\overline{\sigma }(a+vb)=\sigma (a)-vb, \;\;\;\forall \;\;a, b\in B, 
\label{extinvol}
\end{eqnarray}
is an involution for the Cayley-Dickson algebra $\overline{B}(q)$. 
We will show  
that this fact is a particular case of a general result about alternative 
twisted tensor products of algebras, which in turn is analogous to a  
result in \cite{vandaele} about twisted tensor products of 
associative algebras:
\begin{theorem} \label{ext}
Let $A\overline{\otimes }_RB$ be an alternative twisted tensor product of 
algebras, $\sigma _A:A\rightarrow A$ and $\sigma _B:B\rightarrow B$ two 
involutions, and define $\overline{\sigma }:A\otimes B\rightarrow 
A\otimes B$, $\overline{\sigma }:=R\circ (\sigma _B\otimes \sigma _A)\circ 
\tau _{A, B}$. Assume that $R$ satisfies (\ref{braid}) and moreover 
the following conditions hold:
\begin{eqnarray}
&&R(B\otimes A_0)\subseteq A_0\otimes B, \label{inv1}\\
&&\sigma _A(A_0)\subseteq A_0, \label{inv2} \\
&&\overline{\sigma }^2=id_{A\otimes B}. \label{inv3}
\end{eqnarray}
Then $\overline{\sigma }$ is an involution for $A\overline{\otimes }_RB$.
\end{theorem}
\begin{proof}
Note first that $\overline{\sigma }$ is given by the formula 
$\overline{\sigma }(a\otimes b)=\sigma _A(a)_R\otimes \sigma _B(b)_R$, 
for all $a\in A$, $b\in B$, and (\ref{inv3}) together with the fact 
that $\sigma _A$ and $\sigma _B$ are involutions imply 
\begin{eqnarray}
&&\sigma _A(a_R)_r\otimes \sigma _B(b_R)_r=\sigma _A(a)\otimes \sigma _B(b), 
\;\;\;\forall \;\;a\in A,\;b\in B. \label{consinv} 
\end{eqnarray}
In view of (\ref{inv3}), the only thing we need to prove is that 
$\overline{\sigma }$ is antimultiplicative. We will denote by $r$, 
${\cal R}$, $\overline{R}$ some more copies of $R$. First, it is easy to 
see that 
\begin{eqnarray*}
&&\overline{\sigma }((1_A\otimes b)(1_A\otimes b'))=\overline{\sigma }(1_A
\otimes b')\overline{\sigma }(1_A\otimes b), \;\;\;\forall \;\;b, b'\in B.
\end{eqnarray*} 
Take now $a'\in A_0$, $b, b'\in B$. By using (\ref{inv1}) and (\ref{inv2}), 
we compute:
\begin{eqnarray*}
\overline{\sigma }((1_A\otimes b)(a'\otimes b'))&=&
\overline{\sigma }(a'_R\otimes b_Rb')\\
&=&\sigma _A(a'_R)_r\otimes \sigma _B(b_Rb')_r\\
&=&\sigma _A(a'_R)_r\otimes (\sigma _B(b')\sigma _B(b_R))_r\\
&\overset{(\ref{atm3})}{=}&(\sigma _A(a'_R)_{{\cal R}})_r\otimes 
\sigma _B(b_R)_{{\cal R}}\sigma _B(b')_r\\
&\overset{(\ref{consinv})}{=}&\sigma _A(a')_r\otimes 
\sigma _B(b)\sigma _B(b')_r\\
&=&(\sigma _A(a')_r\otimes \sigma _B(b')_r)(1_A\otimes \sigma _B(b))\\
&=&\overline{\sigma }(a'\otimes b')\overline{\sigma }(1_A\otimes b).
\end{eqnarray*}
Take now $a\in A_0$, $b, b'\in B$; by using (\ref{inv2}), we compute: 
\begin{eqnarray*}
\overline{\sigma }((a\otimes b)(1_A\otimes b'))&=&
\overline{\sigma }(a\otimes b'b)\\
&=&\sigma _A(a)_R\otimes \sigma _B(b'b)_R\\
&=&\sigma _A(a)_R\otimes (\sigma _B(b)\sigma _B(b'))_R\\
&\overset{(\ref{atm3})}{=}&(\sigma _A(a)_R)_r\otimes \sigma _B(b')_R
\sigma _B(b)_r\\
&\overset{(\ref{braid})}{=}&(\sigma _A(a)_R)_r\otimes \sigma _B(b')_r
\sigma _B(b)_R\\
&=&(1_A\otimes \sigma _B(b'))(\sigma _A(a)_R\otimes \sigma _B(b)_R)\\
&=&\overline{\sigma }(1_A\otimes b')\overline{\sigma }(a\otimes b).
\end{eqnarray*}
Finally, take $a, a'\in A_0$, $b, b'\in B$. Again by using (\ref{inv1}) and 
(\ref{inv2}) we compute: 
\begin{eqnarray*}
\overline{\sigma }((a\otimes b)(a'\otimes b'))&=&
\overline{\sigma }(aa'_R\otimes b'b_R)\\
&=&\sigma _A(aa'_R)_r\otimes \sigma _B(b'b_R)_r\\
&=&(\sigma _A(a'_R)\sigma _A(a))_r\otimes (\sigma _B(b_R)\sigma _B(b'))_r\\
&\overset{(\ref{atm2})}{=}&\sigma _A(a'_R)_{{\cal R}}\sigma _A(a)_r
\otimes ((\sigma _B(b_R)\sigma _B(b'))_{{\cal R}})_r\\
&\overset{(\ref{atm3})}{=}&(\sigma _A(a'_R)_{\overline{R}})_{{\cal R}}
\sigma _A(a)_r\otimes (\sigma _B(b')_{\overline{R}}
\sigma _B(b_R)_{{\cal R}})_r\\
&\overset{(\ref{braid})}{=}&
(\sigma _A(a'_R)_{\overline{R}})_{{\cal R}}
\sigma _A(a)_r\otimes (\sigma _B(b')_{{\cal R}}
\sigma _B(b_R)_{\overline{R}})_r\\
&\overset{(\ref{consinv})}{=}&
\sigma _A(a')_{{\cal R}}
\sigma _A(a)_r\otimes (\sigma _B(b')_{{\cal R}}
\sigma _B(b))_r\\
&\overset{(\ref{atm3})}{=}&
\sigma _A(a')_{{\cal R}}(\sigma _A(a)_R)_r\otimes \sigma _B(b)_R
(\sigma _B(b')_{{\cal R}})_r\\
&=&(\sigma _A(a')_{{\cal R}}\otimes \sigma _B(b')_{{\cal R}})
(\sigma _A(a)_R\otimes \sigma _B(b)_R)\\
&=&\overline{\sigma }(a'\otimes b')\overline{\sigma }(a\otimes b),
\end{eqnarray*}
finishing the proof.
\end{proof}

\begin{remark}
There exists at least one natural case when (\ref{inv3}) holds automatically. 
Namely, let $A$, $B$ be algebras with $B$ commutative; then obviously 
the flip map $\tau _{A, B}$ is an alternative twisting map. If 
$\sigma _A:A\rightarrow A$ and $\sigma _B:B\rightarrow B$ are two 
involutions, then obviously $\overline{\sigma }$ is given by 
$\overline{\sigma }(a\otimes b)=\sigma _A(a)\otimes \sigma _B(b)$, for all 
$a\in A$, $b\in B$, and it clearly satisfies (\ref{inv3}).   
\end{remark}

Consider now the description $\overline{B}(q)\simeq  
C(K, q)\overline{\otimes }_RB$ of the Cayley-Dickson process as an 
alternative twisted tensor product. We denote $\sigma _B=\sigma $, 
$A=C(K, q)$, and consider the involution $\sigma _A$ of $A$ given by 
$\sigma _A(1)=1$, $\sigma _A(v)=-v$. Then one can easily check that the 
hypotheses of Theorem \ref{ext} are fulfilled, and the involution 
$\overline{\sigma }$ on $A\overline{\otimes }_RB$ provided by the Theorem 
coincides with the involution (\ref{extinvol}) via the identification   
$\overline{B}(q)\simeq C(K, q)\overline{\otimes }_RB$.

We recall from \cite{AM2} that if $A$ is an algebra and 
$\sigma :A\rightarrow A$ is an involutive automorphism, then $\sigma $ 
can be lifted to an involutive automorphism $\overline{\sigma }:Cl(A)
\rightarrow Cl(A)$, $\overline{\sigma }(a+bv)=\sigma (a)-\sigma (b)v$, 
for all $a, b\in A$. It is easy to see that this fact is a particular 
case of the following result about alternative twisted tensor 
products, which in turn is the analogue of a result in \cite{borowiec} 
for twisted tensor products of associative algebras: 
\begin{proposition}
Let $A\overline{\otimes }_RB$ and $E\overline{\otimes }_TF$ be two 
alternative twisted tensor products of algebras and $f:A\rightarrow E$,  
$g:B\rightarrow F$ two algebra maps such that $f(A_0)\subseteq E_0$ and 
$(f\otimes g)\circ R=T\circ (g\otimes f)$. Then 
$f\otimes g:A\overline{\otimes }_RB\rightarrow 
E\overline{\otimes }_TF$ is an algebra map.  
\end{proposition}   
\begin{proof}
A straightforward computation.
\end{proof}
%%%%%%%%%%%%%%%%%%%%%%%%%%%%%%%%%%%%%%%%%%%
\section{A class of examples}
%%%%%%%%%%%%%%%%%%%%%%%%%%%%%%%%%%%%%%%%%%%%
\setcounter{equation}{0}
%%%%%%%%%%%%%%%%%%%%%%%%%%%%%%%%%%%%%%%%%%%%%%%%%
${\;\;\;\;}$
We recall that an algebra $D$ is called {\em alternative} if $(xx)y=x(xy)$ 
and $x(yy)=(xy)y$, for all $x, y\in D$. These identities are called 
the left and respectively right alternative laws. 
We also recall that an involution $\sigma $ on an algebra $B$ is called 
{\em strong} if $b+\sigma (b)\in K\cdot 1_B$ and 
$b\sigma (b)\in K\cdot 1_B$ for all $b\in B$. In this case we denote as 
usual $b+\sigma (b)=t(b)1_B$ and $b\sigma (b)(=\sigma (b)b)=n(b)1_B$. 
The maps $t, n:B\rightarrow K$ are called the {\em trace} and respectively 
{\em norm} of $B$.    
\begin{proposition}
(i) Let $A=K\cdot 1_A\oplus A_0$ be an algebra satisfying the condition
\begin{eqnarray}
&&A_0\cdot A_0\subseteq K\cdot 1_A, \label{grad}
\end{eqnarray}
let $B$ be an algebra and $\sigma :B\rightarrow B$ an involution. Then the 
map $R:B\otimes A\rightarrow A\otimes B$ defined by 
\begin{eqnarray*}
&&R(b\otimes 1_A)=1_A\otimes b, \;\;\;R(b\otimes a)=a\otimes \sigma (b), 
\;\;\;\forall \;\;a\in A_0,\; b\in B, 
\end{eqnarray*}
is an alternative twisting map.\\
(ii) If moreover $A$ and $B$ are alternative and $\sigma $ is strong, 
then, for all $b, b'\in B$ and all ``homogeneous'' $a, a'\in A$  
(i.e. $a$ and $a'$ belong either to $K\cdot 1_A$ or to $A_0$), the 
left and right alternative laws for tensor monomials hold: 
\begin{eqnarray}
&&[(a\otimes b)(a\otimes b)](a'\otimes b')=(a\otimes b)[(a\otimes b)
(a'\otimes b')], \label{alternative1}\\
&&(a\otimes b)[(a'\otimes b')(a'\otimes b')]=[(a\otimes b)(a'\otimes b')]
(a'\otimes b').\label{alternative2}
\end{eqnarray}
\end{proposition}
\begin{proof}
(i) Follows by a direct computation; note that (\ref{grad})  
is used for proving (\ref{atm2}).\\
(ii) We check the left alternative law, while the proof of the right 
alternative law is similar and left to the reader. 
Note first that the alternativity of $B$ implies 
\begin{eqnarray}
&&(b'\sigma (b))b=n(b)b', \label{consalt1} \\
&&(b'b)\sigma (b)=n(b)b', \label{consalt2}
\end{eqnarray}
for all $b, b'\in B$. Indeed, we have
\begin{eqnarray*}
(b'\sigma (b))b&=&(b'(t(b)1_B-b))b\\
&=&t(b)b'b-(b'b)b\\
&=&t(b)b'b-b'(bb)\\
&=&b'(t(b)b-bb)\\
&=&b'(\sigma (b)b)\\
&=&n(b)b',
\end{eqnarray*}
and similarly for (\ref{consalt2}). 

We have to prove that (\ref{alternative1}) holds for all homogeneous 
$a, a'\in A$ and all $b, b'\in B$. This is obvious if $a=a'=1_A$, so we 
only have three cases to analyze:\\
case 1: $a, a'\in A_0$; we compute: 
\begin{eqnarray*}
[(a\otimes b)(a\otimes b)](a'\otimes b')&=&
(a^2\otimes b\sigma (b))(a'\otimes b')\\
&=&(a^2\otimes n(b)1_B)(a'\otimes b')\\
&\overset{(\ref{grad})}{=}&a^2a'\otimes n(b)b', 
\end{eqnarray*}  
\begin{eqnarray*}
(a\otimes b)[(a\otimes b)(a'\otimes b')]&=&(a\otimes b)(aa'\otimes 
b'\sigma (b))\\
&\overset{(\ref{grad})}{=}&a(aa')\otimes (b'\sigma (b))b\\
&\overset{(\ref{consalt1})}{=}&a^2a'\otimes n(b)b'.
\end{eqnarray*}
case 2: $a\in A_0$, $a'=1_A$; we compute:  
\begin{eqnarray*}
[(a\otimes b)(a\otimes b)](1_A\otimes b')&=&
(a^2\otimes n(b)1_B)(1_A\otimes b')\\
&=&a^2\otimes n(b)b', 
\end{eqnarray*}  
\begin{eqnarray*}
(a\otimes b)[(a\otimes b)(1_A\otimes b')]&=&(a\otimes b)(a\otimes  
b'b)\\
&=&a^2\otimes (b'b)\sigma (b)\\
&\overset{(\ref{consalt2})}{=}&a^2\otimes n(b)b'.
\end{eqnarray*}
case 3: $a=1_A$, $a'\in A_0$; we compute: 
\begin{eqnarray*}
[(1_A\otimes b)(1_A\otimes b)](a'\otimes b')&=&
(1_A\otimes b^2)(a'\otimes b')\\
&=&a'\otimes \sigma (b^2)b'\\
&=&a'\otimes \sigma (b)^2b', 
\end{eqnarray*}  
\begin{eqnarray*}
(1_A\otimes b)[(1_A\otimes b)(a'\otimes b')]&=&(1_A\otimes b)(a'\otimes  
\sigma (b)b')\\
&=&a'\otimes \sigma (b)(\sigma (b)b')\\
&=&a'\otimes \sigma (b)^2b',
\end{eqnarray*}
finishing the proof.
\end{proof}
\begin{remark}
Although the alternative laws for homogeneous tensor monomials hold, 
in general this  
does {\bf not} imply that $A\overline{\otimes }_RB$ is an alternative 
algebra. A concrete counterexample is provided by the algebra of 
sedenions, which is not alternative although is an alternative twisted 
tensor product between $C(K, q)$ and the (alternative) algebra 
of octonions.  
\end{remark}

Let now $B$ be an algebra endowed with a strong involution  
$\sigma :B\rightarrow B$, let $q, r\in K\backslash \{0\}$ 
and consider the 3-dimensional algebra $A$ over $K$ with basis 
$\{1, v, z\}$, where $1$ is the unit and    
$v^2=q1$, $z^2=r1$, $vz=zv=0$. Note that $A$ is {\bf not} associative, 
not even alternative (as $qz=(vv)z\neq v(vz)=0$).  
If we define $A_0=Kv\oplus Kz$, then  
obviously we have $A_0\cdot A_0\subseteq K\cdot 1_A$, so we can define 
the alternative twisting map $R:B\otimes A\rightarrow A\otimes B$, 
$R(b\otimes 1_A)=1_A\otimes b$, $R(b\otimes v)=v\otimes \sigma (b)$,  
$R(b\otimes z)=z\otimes \sigma (b)$, for all $b\in B$. We denote the algebra 
$A\overline{\otimes }_RB$ by $\overline{B}(q, r)$ and we write an 
element of this algebra uniquely as $x=a+vb+zc$, with $a, b, c\in B$. 
Then $dim (\overline{B}(q, r))=3\cdot dim (B)$ and the multiplication in 
$\overline{B}(q, r)$ is given by 
\begin{eqnarray*}
&&(a+vb+zc)(a'+vb'+zc')=(aa'+qb'\sigma (b)+rc'\sigma (c))+ 
v(\sigma (a)b'+a'b)+z(\sigma (a)c'+a'c). 
\end{eqnarray*} 
\begin{proposition}
The map $\overline{\sigma }:\overline{B}(q, r)\rightarrow 
\overline{B}(q, r)$, $\overline{\sigma }(a+vb+zc):=\sigma (a)-vb-zc$, 
for all $a, b, c\in B$, is a strong involution for $\overline{B}(q, r)$. 
Moreover, the trace and norm of the element $x=a+vb+zc\in 
\overline{B}(q, r)$ are given by $t(x)=t(a)$ and $n(x)=n(a)-qn(b)-rn(c)$, 
where $t$ and $n$ are the trace and norm on $B$.   
\end{proposition}
\begin{proof}
The fact that $\overline{\sigma }$ is an involution follows either by a 
direct computation or as a consequence of Theorem \ref{ext}, and the fact 
that $\overline{\sigma }$ is strong, together with the explicit formulae for 
the trace and the norm, follow easily by direct computation.   
\end{proof}
\begin{remark}
Obviously the Cayley-Dickson algebras $\overline{B}(q)$ and 
$\overline{B}(r)$ are subalgebras in $\overline{B}(q, r)$. On the other hand, 
if we consider the algebra $\overline{\overline{B}(q)}(r)$ obtained by 
applying the Cayley-Dickson process to $\overline{B}(q)$, it might be 
tempting to believe that $\overline{B}(q, r)$ is a subalgebra in 
$\overline{\overline{B}(q)}(r)$, but this is {\bf not} true, since 
$vz=0$ in $\overline{B}(q, r)$ while $vz\neq 0$ in 
$\overline{\overline{B}(q)}(r)$.  
\end{remark}

Let $B$ be an associative algebra and $\sigma :B\rightarrow B$ a strong 
involution. It is well-known that in this case the Cayley-Dickson algebra 
$\overline{B}(q)$ is alternative. On the other hand, since the 
3-dimensional algebra $A$ defined above is a subalgebra in 
$\overline{B}(q, r)$ and $A$ is not alternative, it follows that 
$\overline{B}(q, r)$ is {\em never} alternative. However, since 
$\overline{B}(q, r)$ is endowed with a strong involution (hence is of 
degree two), it follows that it is always power-associative. 

We will see that the algebras  
$\overline{B}(q, r)$ share some more common properties with the 
Cayley-Dickson algebras $\overline{B}(q)$. We recall that an algebra $D$ 
is called {\em flexible} if $(xy)x=x(yx)$ for all $x, y\in D$. By 
\cite{schaf} we know that all Cayley-Dickson algebras are flexible. We will 
prove that a similar result holds for the algebras $\overline{B}(q, r)$. 
We need to recall some formulae from \cite{schaf}, which are valid for any 
flexible algebra $B$ endowed with a strong involution $\sigma $: 
\begin{eqnarray}
&&(xy)\sigma (y)=\sigma (y)(yx)=y(\sigma (y)x)=(x\sigma (y))y, 
\label{sch1} \\
&&(u\sigma (y))x+y(\sigma (u)x)=x(y\sigma (u))+(xu)\sigma (y), \label{sch2} 
\end{eqnarray}
for all $x, y, u\in B$. Also, by writing $\sigma (x)=t(x)1_B-x$, one can 
easily check the following formula:
\begin{eqnarray}
&&(\sigma (u)\sigma (x))y+x(uy)=\sigma (x)(\sigma (u)y)+(ux)y. \label{sch3}
\end{eqnarray} 
\begin{proposition}
Let $B$ be an algebra and $\sigma :B\rightarrow B$  
a strong involution. Then $\overline{B}(q, r)$ is flexible  
if and only if $B$ is flexible.    
\end{proposition}
\begin{proof}
If $\overline{B}(q, r)$ is flexible then obviously $B$ is  
flexible, as it is a subalgebra of $\overline{B}(q, r)$. 
Conversely, assume that $B$ is flexible and let 
$X=a+vb+zc$ and $Y=a'+vb'+zc'$ be two elements in $\overline{B}(q, r)$. 
We can easily compute: 
\begin{eqnarray*}
(XY)X&=&(aa')a+q((b'\sigma (b))a+b(\sigma (b')a))+r((c'\sigma (c))a
+c(\sigma (c')a))\\
&&+qb(\sigma (b)\sigma (a'))+rc(\sigma (c)\sigma (a'))\\
&&+v((\sigma (a')\sigma (a))b+a(a'b)+q(b\sigma (b'))b+r(c\sigma (c'))b+
a(\sigma (a)b'))\\
&&+z((\sigma (a')\sigma (a))c+a(a'c)+q(b\sigma (b'))c+r(c\sigma (c'))c+
a(\sigma (a)c')), 
\end{eqnarray*}
\begin{eqnarray*}
X(YX)&=&a(a'a)+q(a(b\sigma (b'))+(ab')\sigma (b))+r(a(c\sigma (c'))+
(ac')\sigma (c))\\
&&+q(\sigma (a')b)\sigma (b)+r(\sigma (a')c)\sigma (c)\\
&&+v(\sigma (a)(\sigma (a')b)+(a'a)b+q(b\sigma (b'))b+r(c\sigma (c'))b+
\sigma (a)(ab'))\\
&&+z(\sigma (a)(\sigma (a')c)+(a'a)c+q(b\sigma (b'))c+r(c\sigma (c'))c+
\sigma (a)(ac')),
\end{eqnarray*}
and we immediately obtain $(XY)X=X(YX)$ by using the flexibility of $B$ 
together with (\ref{sch1}), (\ref{sch2}), (\ref{sch3}). 
\end{proof}
\begin{remark}
In particular, if we take $B$ to be an associative or alternative algebra 
with a strong involution, we obtain an  
infinite family of flexible algebras, of dimensions 
$3^n\cdot dim (B)$,  
for all positive integers $n$.  
\end{remark}

Assume now $char (K)\neq 2$ and let $B$ be an algebra endowed 
with a strong involution $\sigma $. We recall that the norm $n$ on $B$ 
is called {\em nondegenerate} if the associated symmetric bilinear 
form 
\begin{eqnarray*}
&&(x, y)=\frac{1}{2}(n(x+y)-n(x)-n(y)), \;\;\;\forall \;\;x, y\in B, 
\end{eqnarray*} 
is nondegenerate (i.e. if $(x, y)=0$ for all $y\in B$ then $x=0$). It is 
well-known (see \cite{schafer}, p. 70) that all Cayley-Dickson algebras have 
nondegenerate norms. We prove a similar result for the algebras 
$\overline{B}(q, r)$: 
\begin{proposition}
If $B$ has nondegenerate norm, then $\overline{B}(q, r)$ has also 
nondegenerate norm.
\end{proposition}
\begin{proof}
Let $x=a+vb+zc\in \overline{B}(q, r)$ such that $(x, y)=0$ for all 
$y=a'+vb'+zc'\in \overline{B}(q, r)$. It is easy to see that 
\begin{eqnarray*}
(x, y)&=&\frac{1}{2}(n(x+y)-n(x)-n(y))\\
&=&\frac{1}{2}(n(a+a')-qn(b+b')-rn(c+c')-n(a)+qn(b)+rn(c)\\
&&-n(a')+qn(b')+rn(c'))\\
&=&(a, a')-q(b, b')-r(c, c'), 
\end{eqnarray*} 
thus $(x, y)=0$ for all $y$ implies $(a, a')=q(b, b')+r(c, c')$ for all 
$a', b', c'\in B$. By taking $b'=c'=0$ we obtain $(a, a')=0$ for all 
$a'\in B$, and since the norm on $B$ is nondegenerate this implies 
$a=0$. Similarly, by taking $a'=c'=0$ and then $a'=b'=0$ and using 
the fact that $q\neq 0$, $r\neq 0$ and again the nondegeneracy of the norm 
on $B$, we obtain $b=0$ and respectively $c=0$; that is, $x=0$. 
\end{proof}
%%%%%%%%%%%%%%%%%%%%%%%%%%%%%%%%%%%%%%%%%%%%%%%%%%%%%%%%%%%%%%%%%%%%%%%%%
\begin{center}
ACKNOWLEDGEMENTS
\end{center}
We would like to thank the referee for some useful comments and suggestions.
%%%%%%%%%%%%%%%%%%%%%%%%%%%%%%%%%%%%%%%%%%%%%%%%%%%%%%%%%%%%%%%%%%%%%%%%
%%%%%%%%%%%%%%%%%%%%%%%%%%%%%%%%%%%%%%%%%%%%%%%%%%%%%%%

\end{document}